\documentclass[11pt]{amsart}
\usepackage{amsmath,amssymb,amsthm}
\usepackage[margin=1in]{geometry}
\usepackage[hidelinks]{hyperref}
\newtheorem{theorem}{Theorem}[section]
\newtheorem{lemma}[theorem]{Lemma}
\newtheorem{corollary}[theorem]{Corollary}
\newtheorem{proposition}[theorem]{Proposition}
\theoremstyle{remark}
\newtheorem{remark}[theorem]{Remark}
\newtheorem{example}[theorem]{Example}
\newtheorem{question}[theorem]{Question}
\newcommand{\G}{\mathcal G}
\title[Examples to the non-abelian Gyarmati--Matolcsi--Ruzsa problem]{Why the non-abelian 
Gyarmati--Matolcsi--Ruzsa inequality does not generally work:\\ mechanisms and explicit examples}
\author{Christian Elsholtz}
\address{Graz University of Technology, Institut f\"ur Analysis und Zahlentheorie, Kopernikusgasse 24, 8010 Graz, Austria}
\date{\today}
\begin{document}
\begin{abstract}
For finite sets in a commutative semigroup, Gyarmati, Matolcsi and Ruzsa proved
\[ |A_1+\dots+A_n|^{n-1}\ \le\ \prod_{i=1}^n|A_1+\dots+A_{i-1}+A_{i+1}+\dots+A_n| . \]
They asked whether in an arbitrary group
\[ |A_1\cdots A_n|^{n-1}\ \le\ \prod_{i=1}^n\ \max_{a\in A_i}\,|A_1\cdots A_{i-1}\,a\,A_{i+1}\cdots A_n| , \]
where each omitted set is replaced by the largest slice obtained by fixing one of its elements. 
Nair, Namayandeh and Zhu proved that this inequality does not hold when 
$n=4$ and $G=D_4^{14}$, of group order $8^{14}$. 
We develop its double-coset mechanism and construct counterexamples in many other groups.

For $n=4$, the smallest group found has order $518\,400$, namely $A_6^2\times C_2^2$; no absolute minimality is claimed. For $n=6$ there is a witness in a group of order $1728$, for $n=10$ a class-three $3$-group of order $243$ suffices, and for $n=4$ a single class-two $2$-group of order $2^{23}$ works. We also obtain torsion-free examples. The violation ratio $R_n=|A_1\cdots A_n|^{n-1}/\prod_iN_i$ is greater than one exactly when the inequality fails; explicit four-set constructions attain $R_4\approx1.29$ in $A_8^2$, and a larger ratio of about $1.43$ is reported for $A_{10}^2$.

For every $n\ge4$, the inequality holds in all direct powers of $G$ precisely when $G$ is abelian. For $n\ge2$ the optimal exponent $\theta_n$ is the infimum of $\log(\prod_iN_i)/\log|A_1\cdots A_n|$ over products of size at least two. We obtain $2\le\theta_n\le n-1.0768$, and an infinite family in $C_3*C_3$ shows $\theta_n\le n-1-\delta$ for every $\delta<5-2\log_25\approx0.3561$ and all sufficiently large $n$.
\end{abstract}
\maketitle

\section{Introduction}

\subsection{The theorem of Gyarmati, Matolcsi and Ruzsa}
Let $A_1,\dots,A_n$ be finite nonempty sets in a commutative semigroup, $S=A_1+\dots+A_n$, and let $S_i$ be the same sum with $A_i$ omitted. Gyarmati, Matolcsi and Ruzsa [GMR, Theorem~1.2] proved
\begin{equation}\label{eq:gmr}
|S|^{n-1}\le\prod_{i=1}^n|S_i| ;
\end{equation}
the case $n=3$ is due to Ruzsa [R07, GMR]. Inequality \eqref{eq:gmr} is an additive analogue of the Loomis--Whitney projection inequality [LW] and of Han's entropy inequality [Han]; equality holds for sets on independent coordinate axes of $\mathbb Z^n$. Related sumset and entropy inequalities are discussed in [BB, MMT, R, T, LN].

\subsection{The non-commutative question}
In a non-commutative group $(\G,\circ)$, replacing a summand by a single element $a$ no longer gives a translate of a fixed set. For an inner index $i$, the set $A_1\cdots A_{i-1}\,a\,A_{i+1}\cdots A_n$ genuinely depends on $a$. [GMR] therefore replace $|S_i|$ by the largest such slice,
\[ N_i=\max_{a\in A_i}\bigl|A_1\cdots A_{i-1}\,a\,A_{i+1}\cdots A_n\bigr| \]
(for $i=1$ and $i=n$ this is simply $|A_2\cdots A_n|$, respectively $|A_1\cdots A_{n-1}|$). They asked [GMR, Problem~1.3] whether
\begin{equation}\label{eq:nagmr}
|A_1\cdots A_n|^{n-1}\ \le\ \prod_{i=1}^nN_i
\end{equation}
holds in every group. For $n=3$ the non-abelian inequality is also valid: Ruzsa's three-set argument can be adapted to groups (see [R07, GMR, NNZ]). The adaptation uses inverses and does not by itself give a statement for general non-commutative semigroups.

\smallskip\noindent\emph{Reasons to believe \eqref{eq:nagmr}.} Its abelian and geometric analogues (Han, Loomis--Whitney, Shearer) are purely combinatorial and make no use of commutativity. The case $n=3$ holds in every group. In the abelian case $N_i=|S_i|$, so \eqref{eq:nagmr} is the most natural non-commutative formulation.

\smallskip\noindent\emph{Reasons to doubt it.} Each $N_i$ is a maximum over \emph{single} elements of $A_i$, so it cannot see effects that need two summands to cooperate. In a non-abelian group such effects are forced by double cosets. Let $H_1,H_2$ be subgroups, $y\in H_2$ and $z\in H_1$. Then $H_1yH_2=H_1zH_2=H_1H_2$: each of $y$ and $z$ alone is harmless. But the product $H_1\,yz\,H_2$ can be a new double coset disjoint from $H_1H_2$. In an abelian group this cannot happen, because $H_1+H_2$ is a subgroup and $H_1+y+z+H_2=H_1+H_2$. Taking $A_1=H_1$, $A_4=H_2$, and letting $A_2,A_3$ contain $y$, $z$ (sparsely, so that the joint event is rare) makes $|S|$ larger than every slice ``predicts''. This is the mechanism behind every counterexample below.

\subsection{The breakthrough of Nair, Namayandeh and Zhu}\label{sec:nnz}
A counterexample consists of sets with $|S|^{n-1}>\prod_iN_i$. It is convenient to measure the failure by the \emph{violation ratio}
\[ R_n=R_n(A_1,\dots,A_n)=\frac{|A_1\cdots A_n|^{n-1}}{\prod_{i=1}^nN_i}, \]
so that \eqref{eq:nagmr} fails exactly when $R_n>1$; for $n=4$ we write $R=R_4$. 
The pioneering work [NNZ] resolved the conjecture by exhibiting the first explicit failure for $n=4$. 
Finding this first witness was a major conceptual step; verifying its cardinalities is elementary. 
The method of types [NNZ, Proposition~B.1] relates the problem to a maximal-coupling entropy 
inequality: for fixed $n$, \eqref{eq:nagmr} holds for all finite sets in all groups (it suffices 
to consider all direct powers $\G^m$ of a fixed group) if and only if
\[ (n-1)\max_{q\in\Pi}H_q(X_1\circ\dots\circ X_n)\ \le\ 
\sum_{i=1}^n\max_{q\in\Pi}H_q(X_1\circ\dots\circ X_n\mid X_i) \]
holds for all finitely supported marginals. Here $\Pi$ is the set of couplings with the given 
marginals, and each maximum is taken separately. With AI assistance they found an entropic counterexample 
in the dihedral group $D_4$ of order $8$: two binary variables are uniform and two are $\varepsilon$-biased, 
and the gap is $h_2(\varepsilon)-4\varepsilon>0$ for $0<\varepsilon<1/17$. They then built explicit 
sets in $D_4^m$ with
\[ \frac{|S|^3}{N_1N_2N_3N_4}=\frac{(m+1)^3}{4(2m+1)^2}, \]
which exceeds $1$ for $m=14$ ($3375/3364$), in a group of order $8^{14}\approx4.4\cdot10^{12}$. 
Although AI tools helped to \emph{find} it, the final witness is short and its cardinalities can be checked 
by hand, by classifying the patterns of the single exceptional coordinate.

\subsection{Main results}
\begin{enumerate}
\item \textbf{Mechanism} (Sections~\ref{sec:mech}--\ref{sec:multi}). Two subgroups $H_1,H_2$ and rare elements $z\in H_1$, $y\in H_2$ with $yz\notin H_1H_2$ (a \emph{double-coset collapse}) give explicit ratios (Theorems~\ref{thm:binary}, \ref{thm:multi}); NNZ is an instance, and $D_4^{13}$ already works.
\item \textbf{All non-abelian groups} (Section~\ref{sec:all}). Every non-abelian $G$ fails the four-set inequality in $G^{42}$ (in $G^{39}$ if it has non-commuting elements of finite order). So the inequality holds in all direct powers of $G$ iff $G$ is abelian.
\item \textbf{Four sets} (Section~\ref{sec:witness}). Explicit witnesses in $A_6^2\times C_2^2$ (order $518\,400$), $A_8^2$ (ratio $1.29$), $A_{10}^2$ (ratio $1.43$), and torsion-free ones in $F_2\times F_2$.
\item \textbf{More sets, smaller groups} (Section~\ref{sec:n}). A five-double-coset criterion (Theorem~\ref{thm:crit}) gives six-set counterexamples in $A_5^2$, in $\mathrm{PSL}(2,3^f)^2$ for all $f\ge2$, and in a central quotient of order $1728$; splitting the outer sets gives $7$- and $8$-set examples with no singleton sets and ratio $1.1675$. A ten-set counterexample lives in the free class-$3$ quotient of $C_3*C_3$, of order $243$ (Theorem~\ref{thm:243}).
\item \textbf{An infinite family and the exponent} (Section~\ref{sec:fam}). For $n=4s+2$ an explicit configuration in $C_3*C_3$ (realised in $\mathrm{PSL}(2,p)$) has ratio $\sim\frac13(1.6384)^s$ (Theorem~\ref{thm:fam}). Hence $(n-1)-\theta_n\ge0.11$ for $n\ge18$ and $\liminf_n((n-1)-\theta_n)\ge5-2\log_25\approx0.356$, where $\theta_n$ is the optimal exponent.
\end{enumerate}
The results below follow from the stated arguments or from exact finite certificates for the specific instances identified; several existence statements in Sections~\ref{sec:witness}--\ref{sec:n} use explicit computations. Unsuccessful restricted searches are computational evidence, not proofs of group-wide non-existence (Section~\ref{sec:method}).

\subsection{Notation and conventions}
Throughout, $e$ is the identity, $\G$ a group, $A_i$ finite nonempty sets, $S=A_1\cdots A_n$, and $N_i$ the largest product obtained by fixing one element of $A_i$. $R_n=|S|^{n-1}/\prod_iN_i$ is the violation ratio of Section~\ref{sec:nnz}, and $\theta_n$ the optimal exponent (defined in Section~\ref{sec:exp}). For $n=4$ we write $R=R_4$. $C_m$ is the cyclic group of order $m$, $D_m$ the dihedral group of order $2m$, $S_m$ and $A_m$ the symmetric and alternating groups, $Q_8$ the quaternion group; $\mathrm{SL}$, $\mathrm{PSL}$ are the special and projective special linear groups, $\mathrm{Sz}(q)$ a Suzuki group, $F_2$ the free group on two generators, and $C_3*C_3=\langle a,b\mid a^3=b^3=e\rangle$ the free product of two cyclic groups of order $3$, an infinite group in which every element has a unique reduced form as an alternating word in the non-identity powers of $a$ and of $b$. $H^\#=H\setminus\{e\}$, $H'$ is the derived subgroup, $Z(H)$ the centre, $\gamma_i(H)$ the lower central series, and $HxJ$ a double coset. $\G^m$ is a direct power, and $x^{(i)}$ has $x$ in coordinate $i$ and $e$ elsewhere. For subgroups $H_1,H_2$ of orders $k,l$, the quantity $\rho=|H_1yzH_2|/(kl)$ is introduced in Theorem~\ref{thm:binary}.

\section{Preliminaries}\label{sec:prelim}

Let $\G$ be a group and $A_1,\dots,A_n\subseteq\G$ finite and nonempty, $S=A_1\cdots A_n$, and
\[ N_i=\max_{a\in A_i}|A_1\cdots A_{i-1}aA_{i+1}\cdots A_n|. \]
For $n=4$ write $(A,B,C,D)$, so $N_1=|BCD|$, $N_4=|ABC|$, $N_2=\max_b|AbCD|$ and $N_3=\max_c|ABcD|$. The \emph{ratio} is $R=R_\G(A,B,C,D)=|S|^3/(N_1N_2N_3N_4)$; a counterexample is $R>1$.

\begin{remark}[Normalised form]
With $a=|A|$, $b=|B|$, $c=|C|$, $d=|D|$ put
\[ t=\frac{|ABCD|}{abcd},\quad t_A=\frac{N_1}{bcd},\quad t_B=\frac{N_2}{acd},\quad t_C=\frac{N_3}{abd},\quad t_D=\frac{N_4}{abc}. \]
Then exactly $R=t^3/(t_At_Bt_Ct_D)$. Thus a counterexample needs a specific imbalance: the full product must retain relatively more distinct outputs than the four slice products. Many collisions alone are not enough.
\end{remark}

\begin{lemma}[Padding]\label{lem:pad}
A counterexample for $n=4$ gives one for every $n\ge4$ in the same group: add $A_5=\dots=A_n=\{e\}$.
\end{lemma}
\begin{proof}
$S$ and $N_1,\dots,N_4$ are unchanged and $N_5=\dots=N_n=|S|$, so the inequality becomes $|S|^3\le N_1N_2N_3N_4$.
\end{proof}

\begin{lemma}[Translations]\label{lem:transl}
Replacing $(A_i,A_{i+1})$ by $(A_ih,h^{-1}A_{i+1})$, $A_1$ by $gA_1$, or $A_n$ by $A_ng$ leaves $|S|$ and all $N_i$ unchanged. Hence one may assume $e\in A_i$ for all $i$.
\end{lemma}

Notation: in $G^m$, $x^{(i)}$ denotes the element with $x$ in coordinate $i$ and $e$ elsewhere.

\begin{proposition}[Subgroup inheritance]\label{prop:sub}
If $H\le\G$ contains a counterexample, then so does $\G$, with the same ratio.
\end{proposition}
\begin{proof}
The products and all slices of subsets of $H$ are the same sets in $\G$.
\end{proof}

\begin{proposition}[Multiplicativity]\label{prop:mult}
If $(A_i,B_i,C_i,D_i)$ are configurations in $\G_i$ ($i=1,2$), then the coordinatewise products in $\G_1\times\G_2$ have ratio $R_1R_2$.
\end{proposition}
\begin{proof}
The full product is the product of the two full products. For each slice, the maximum over pairs factorises into the two maxima.
\end{proof}
In particular, $\G\times K$ is a counterexample group for every group $K$ whenever $\G$ is one, and one counterexample with ratio $R$ gives ratio $R^t$ in $\G^t$.

\begin{lemma}[Quotients]\label{lem:quot}
Let $(A_1,\dots,A_n)$ be finite subsets of $\G$, let $N\trianglelefteq\G$, and let $\pi:\G\to\G/N$ be the projection. If $\pi$ is injective on $S=A_1\cdots A_n$, then the ratio for $(\pi(A_1),\dots,\pi(A_n))$ is at least the ratio for $(A_1,\dots,A_n)$.
\end{lemma}
\begin{proof}
$\pi(A_1)\cdots\pi(A_n)=\pi(S)$ has $|S|$ elements, while every slice can only shrink under $\pi$.
\end{proof}
Quotients of counterexample groups need not be counterexample groups in general, because collisions can change $|S|$ and the slices. Injectivity on the full product is exactly what makes the central quotients below valid. They are the only quotients that keep both coordinates:

\begin{lemma}[Twisting quotients]\label{lem:twist}
Let $N\trianglelefteq L\times L$ with $N\cap(L\times\{e\})=N\cap(\{e\}\times L)=\{e\}$. Then $N=\{(z,\varphi(z)):z\in Z_1\}$ for subgroups $Z_1,Z_2\le Z(L)$ and an isomorphism $\varphi:Z_1\to Z_2$.
\end{lemma}
\begin{proof}
For $(x,y)\in N$ and $\ell\in L$, the commutator $[(x,y),(\ell,e)]=([x,\ell],e)$ lies in $N\cap(L\times\{e\})=\{e\}$, so $x\in Z(L)$; similarly $y\in Z(L)$. Both projections of $N$ are injective, so $N$ is the graph of an isomorphism between its projections.
\end{proof}

\section{The mechanism: a double-coset collapse}\label{sec:mech}

\begin{theorem}\label{thm:binary}
Let $H_1,H_2$ be finite subgroups of a group $G$, $|H_1|=k$, $|H_2|=l$, with $H_1\cap H_2=\{e\}$, and let $z\in H_1$, $y\in H_2$ with $yz\notin H_1H_2$. Put $d=|H_1yzH_2|$ and $\rho=d/(kl)\in(0,1]$. For $m\ge2$ define in $G^m$
\[ A=H_1^m,\quad B=\{y^{(i)}:i\le m\},\quad C=\{z^{(j)}:j\le m\},\quad D=H_2^m. \]
Then
\[ |S|=(kl)^{m-1}(kl+md),\ \ N_1=2m\,l^m,\ \ N_2=N_3=(kl)^{m-1}(kl+d),\ \ N_4=2m\,k^m, \]
and hence \[ R=\frac{(1+m\rho)^3}{4m^2(1+\rho)^2}. \]
\end{theorem}

\begin{proof}
Let $K=H_1H_2$; $|K|=kl$ since $H_1\cap H_2=\{e\}$. As $y\in H_2$ and $z\in H_1$, $H_1yH_2=H_1zH_2=K$, while $D_1:=H_1yzH_2$ is a double coset disjoint from $K$. The left cosets $H_2,\,zH_2,\,yzH_2$ are pairwise distinct: $zH_2=H_2$ forces $z\in H_1\cap H_2$, so $yz=y\in K$; $yzH_2=H_2$ forces $yz\in K$; and $zH_2=yzH_2$ means $z^{-1}yz\in H_2$, so $yz=z\cdot z^{-1}yz\in K$. Symmetrically the right cosets $H_1,H_1y,H_1yz$ are distinct.

All products are coordinatewise. For $i\ne j$, $Ay^{(i)}z^{(j)}D=K^m$; for $i=j$ it is $K^{i-1}\times D_1\times K^{m-i}$. Since $m\ge2$, $S$ is the disjoint union of $K^m$ and these $m$ sets, which gives $|S|$. For fixed $i$, $Ay^{(i)}CD=K^m\cup(K^{i-1}\times D_1\times K^{m-i})$, giving $N_2$; $N_3$ is symmetric. Next, $y^{(i)}z^{(j)}D$ equals $H_2^{j-1}\times zH_2\times H_2^{m-j}$ if $i\ne j$ and $H_2^{j-1}\times yzH_2\times H_2^{m-j}$ if $i=j$; these $2m$ sets are pairwise disjoint, so $N_1=|BCD|=2ml^m$ (left multiplication by $a$ is a bijection). $N_4$ is symmetric.
\end{proof}

The equality $\rho=1$ holds exactly when $H_1\cap(yz)H_2(yz)^{-1}=\{e\}$. Then $R=(1+m)^3/(16m^2)>1$ for $m\ge13$ ($2744/2704$ at $m=13$). For $\rho=\frac12,\frac13,\frac14$ the first admissible $m$ are $66$, $183$, $388$.

\begin{corollary}[Involutions]\label{cor:inv}
Let $y,z\in G$ be involutions and let $o$ be the order of $yz$. Take $H_1=\{e,z\}$ and $H_2=\{e,y\}$. If $o\ge4$, then $\rho=1$ and Theorem~\ref{thm:binary} gives a counterexample in $G^{13}$. If $o=3$, then $\rho=\frac12$ and it gives one in $G^{66}$.
\end{corollary}
\begin{proof}
$yz\in\{e,y,z,zy\}$ iff $yz=zy$ iff $o\le2$. Moreover $H_1\cap(yz)H_2(yz)^{-1}\ne\{e\}$ iff $(yz)y(yz)^{-1}=z$ iff $(yz)^3=e$.
\end{proof}

\begin{example}[Simplified Nair--Namayandeh--Zhu construction]
In $D_4=\langle r,s\mid r^4=s^2=e,\ srs=r^{-1}\rangle$ choose $y=s$ and $z=rs$; then $yz=r^{-1}$ has order four. For $m=13$ define four sets in $D_4^{13}$ by
\[ A=\{e,rs\}^{13},\quad B=\{s^{(i)}:1\le i\le13\},\quad C=\{(rs)^{(j)}:1\le j\le13\},\quad D=\{e,s\}^{13}. \]
The counting of Theorem~\ref{thm:binary} gives $|ABCD|=14\cdot4^{13}$, $N_1=N_4=26\cdot2^{13}$ and $N_2=N_3=2\cdot4^{13}$, so $R=2744/2704>1$. More generally, two non-commuting involutions whose product has order at least four give the same construction in the $13$th power of their group; this applies to $D_m$ for $m\ge4$ and to $S_k$ for $k\ge4$. In $S_3$ the product has order three, and Corollary~\ref{cor:inv} requires $66$ coordinates. The original Nair--Namayandeh--Zhu sets use weight at most one rather than weight exactly one, giving $R=(m+1)^3/(4(2m+1)^2)$ and requiring $m=14$.
\end{example}

\begin{remark}[Entropic form]
In the entropy formulation [NNZ, (3)], take $X$ and $W$ uniform on $H_1$ and $H_2$, and $\Pr(Y=y)=\Pr(Z=z)=\varepsilon$, $Y\in\{e,y\}$, $Z\in\{e,z\}$. The same coset facts give the gap
\[ 3\max H(U)-\sum_i\max H(U\mid X_i)\ \ge\ h_2(\varepsilon)-c(\rho)\varepsilon,\qquad c(\rho)=2+2\log_2(1+\rho)+3\log_2(1/\rho), \]
with maximum $\log_2(1+2^{-c(\rho)})$. For $\rho=1$ this is $h_2(\varepsilon)-4\varepsilon$, exactly the NNZ computation. For $S_3$ ($\rho=\frac12$, $c=6.17$) and $A_4$ ($\rho=\frac13$) the bound agrees with a numerical convex optimisation to six decimals.
\end{remark}

\section{Every non-abelian group has a counterexample in a direct power}\label{sec:all}
\begin{theorem}\label{thm:all}
Let $g,h\in G$ be non-commuting elements of finite order. In $G^3$ put $y=(g,e,g)$, $z=(h,h,e)$, $H_1=\langle z\rangle$ and $H_2=\langle y\rangle$. Then the hypotheses of Theorem~\ref{thm:binary} hold with $\rho=1$. Hence $(G^3)^{13}=G^{39}$ contains a counterexample, and so does $G^{39}\times\G'$ for every group $\G'$; by Lemma~\ref{lem:pad} this holds for every $n\ge4$.
\end{theorem}

\begin{proof}
The key invariant: for $w=(w_1,w_2,w_3)$ compare $w_1$ with $w_2w_3$. Elements of $K=H_1H_2$ have the form $z^uy^v=(h^ug^v,h^u,g^v)$, so they satisfy $w_1=w_2w_3$. Elements of $H_1yzH_2$ have the form $z^uyzy^v=(h^ughg^v,\,h^{u+1},\,g^{v+1})$; for these, $w_1=w_2w_3$ iff $gh=hg$. Hence $yz\notin K$. If $z^u=y^v$, the second and third coordinates give $h^u=e=g^v$, so $H_1\cap H_2=\{e\}$. Finally $(u,v)\mapsto z^uyzy^v$ is injective, since coordinates $2$ and $3$ recover $h^u$ and $g^v$. So $|H_1yzH_2|=|H_1||H_2|$, i.e.\ $\rho=1$.
\end{proof}

For elements of infinite order, replace the cyclic subgroups by segments.

\begin{theorem}\label{thm:all2}
Let $g,h\in G$ be arbitrary non-commuting elements and $y,z\in G^3$ as above. Let $P=\langle z\rangle$ if $h$ has finite order, and otherwise $P=\{z^u:0\le u<\kappa\}$. Define $Q\subset\langle y\rangle$ the same way from $g$. In $(G^3)^m$ put $A=P^m$, $B=\{y^{(i)}\}$, $C=\{z^{(j)}\}$, $D=Q^m$. Then
\[ R\ \ge\ \frac{(m+1)^3}{m^2\bigl((1+1/\kappa)^{m-1}+1\bigr)^4}. \]
In particular $m=14$, $\kappa=357$ (or $m=20$, $\kappa=107$) gives a counterexample in $G^{42}$ (resp.\ $G^{60}$).
\end{theorem}

\begin{proof}[Sketch]
Every block entry of an element of $S$ or of a slice has the form $w=z^uy^\beta z^\gamma y^v$ with $\beta,\gamma\in\{0,1\}$, and $w_1=w_2w_3$ unless $\beta=\gamma=1$. So $w$ is determined by its \emph{shadow} $(w_2,w_3,[\beta=\gamma=1])\in\langle h\rangle\times\langle g\rangle\times\{0,1\}$, and the counting of Theorem~\ref{thm:binary} can be repeated on shadows. The only change is that a segment $\bar P$ and its translate $h\bar P$ now overlap in all but one element, which costs the factor $(1+1/\kappa)^{m-1}$ in the slices.
\end{proof}

\begin{corollary}\label{cor:iff}
For a group $G$ the following are equivalent:
(a) for every $n\ge4$ (equivalently, for $n=4$), the inequality $|S|^{n-1}\le\prod N_i$ holds for all finite subsets of all direct powers $G^m$;
(b) $G$ is abelian.
\end{corollary}
\begin{proof}
(b)$\Rightarrow$(a) is the theorem of Gyarmati--Matolcsi--Ruzsa, since $G^m$ is abelian. (a)$\Rightarrow$(b) follows from Theorems~\ref{thm:all} and \ref{thm:all2}.
\end{proof}

\begin{remark}[Property (H)]\label{rem:H}
Say that $G$ has property (H) if it contains cyclic subgroups $P,Q$ with $P\cap Q=\{e\}$ but $PQ\ne QP$. Then $G$ itself (rather than $G^3$) supplies the configuration of Theorem~\ref{thm:binary}, which gives a counterexample in a direct power of $G$. A GAP census found property (H) in $373$ of the $469$ non-abelian groups of order at most $64$; the exceptions include $Q_8$, $C_3\rtimes C_4$, $Q_{16}$, $M_{16}$ and $C_9\rtimes C_3$. This classification was not independently repeated for the present revision.
\end{remark}

\section{The multi-rare construction and its limits}\label{sec:multi}
The binary construction wastes a factor $(|G|/kl)^m$. Allowing several rare elements per coordinate gives much smaller groups.

\begin{theorem}\label{thm:multi}
Let $H_1,H_2\le G$ be finite subgroups with $H_1\cap H_2=\{e\}$, and let $Y'\subseteq H_2\setminus\{e\}$, $Z'\subseteq H_1\setminus\{e\}$ with $|Y'|=r$ and $|Z'|=s$. Assume that the $rs+1$ double cosets $H_1H_2$ and $H_1yzH_2$ ($y\in Y'$, $z\in Z'$) are pairwise distinct and each has size $kl$. In $G^m$ ($m\ge2$) put $A=H_1^m$, $D=H_2^m$, $B=\{y^{(i)}:y\in Y'\}$, $C=\{z^{(j)}:z\in Z'\}$. Then
\[ |S|=(kl)^m(1+mrs),\quad N_2=(kl)^m(1+s),\quad N_3=(kl)^m(1+r), \]
\[ N_1\le ms(r+1)l^m,\qquad N_4\le mr(s+1)k^m, \]
\[ R\ \ge\ \frac{(1+mrs)^3}{m^2\,rs\,(r+1)^2(s+1)^2}. \]
For $m=2$ this exceeds $1$ when, for example, $(r,s)\in\{(5,5),(4,6),(5,6),\dots\}$; for $m=3$ when $r=s=3$.
\end{theorem}
\begin{proof}
Write $K=H_1H_2$ and $D_{yz}=H_1yzH_2$. All products are coordinatewise; $y\in H_2$ and $z\in H_1$ are absorbed: $H_1yH_2=H_1zH_2=K$. For $b=y^{(i)}$ and $c=z^{(j)}$ with $i\ne j$ we get $AbcD=K^m$. For $i=j$ we get $K^{i-1}\times D_{yz}\times K^{m-i}$. Since $m\ge2$, $S$ is the disjoint union of $K^m$ and the $mrs$ sets $K^{i-1}\times D_{yz}\times K^{m-i}$, each of size $(kl)^m$. For fixed $b=y^{(i)}$ the slice $AbCD$ is $K^m\cup\bigcup_{z}K^{i-1}\times D_{yz}\times K^{m-i}$, of size $(kl)^m(1+s)$, and similarly for $N_3$. For $N_1=|BCD|$ note that $y^{(i)}z^{(j)}D$ depends only on $(j,z)$ when $i\ne j$ (since $yH_2=H_2$). This gives at most $ms$ sets of size $l^m$, plus $mrs$ sets for $i=j$, so $N_1\le ms(r+1)l^m$. Symmetrically $N_4\le mr(s+1)k^m$. Combine.
\end{proof}

\begin{remark}[Variants]
Including $e$ in $B$ and $C$ (``weight at most one'') gives the bound $(1+mrs)^3/\bigl((1+ms(r+1))(1+mr(s+1))(1+r)(1+s)\bigr)$; we always use the larger of the two. Different coordinates may carry different rectangles.
\end{remark}

\begin{lemma}[Derived-subgroup obstruction]\label{lem:derived}
In the situation of Theorem~\ref{thm:multi}, let $U=\langle H_1,H_2\rangle$. Then $1+rs\le|U'|$, where $U'$ is the derived subgroup. In particular an $m=2$ configuration needs $|U'|\ge25$, and an $m=3$ configuration needs $|U'|\ge10$.
\end{lemma}
\begin{proof}
The $1+rs$ double cosets are disjoint, lie in $U$ and have size $kl$, so $(1+rs)kl\le|U|$. On the other hand $U/U'$ is abelian and generated by the images of $H_1$ and $H_2$, so $|U/U'|\le|H_1U'/U'|\,|H_2U'/U'|\le kl$. Hence $(1+rs)kl\le|U'|\,kl$.
\end{proof}

\begin{corollary}[Size bound]\label{cor:size}
An $m=2$ configuration needs $(r,s)$ with $(1+2rs)^3>4rs(r+1)^2(s+1)^2$ or the weight-$\le1$ analogue, hence $rs\ge24$ and $|\langle H_1,H_2\rangle|\ge kl(1+rs)\ge875$. The resulting counterexample group, after a central quotient by $Z$ as in Lemma~\ref{lem:quot}, has order at least $875^2/|Z|$, and at least $|U|^2/|Z(U)|$ for $U=\langle H_1,H_2\rangle$. For $m=3$ one needs $kl(1+rs)\ge160$, and for $m=4$ at least $84$.
\end{corollary}

\section{Explicit witnesses for $n=4$}\label{sec:witness}

\subsection{The smallest witness found: $A_6^2\times C_2^2$}

\begin{example}[$A_6\times A_6\times C_2\times C_2$, order $518\,400$]\label{ex:a6}
Let $G_0=A_6\times\langle u\rangle\times\langle v\rangle$ with central involutions $u,v$. Take two ``twisted'' copies of $S_3$, in which the involutions carry a central factor:
\begin{align*}
H_1&=\{e,\ (1\,3\,6)(2\,4\,5),\ (1\,6\,3)(2\,5\,4),\ (3\,6)(4\,5)u,\ (1\,3)(2\,4)u,\ (1\,6)(2\,5)u\},\\
H_2&=\{e,\ (1\,5\,6),\ (1\,6\,5),\ (2\,3)(5\,6)v,\ (1\,5)(2\,3)v,\ (1\,6)(2\,3)v\}.
\end{align*}
Then $H_1\cap H_2=\{e\}$, $\langle H_1,H_2\rangle=G_0$ (order $1440$), and the $26$ double cosets $H_1H_2$ and $H_1yzH_2$ ($y\in H_2^\#$, $z\in H_1^\#$) are distinct and of size $36$. So Theorem~\ref{thm:multi} applies with $r=s=5$, $m=2$, in $G_0\times G_0$. The projection $G_0\times G_0\to A_6\times A_6\times C_2\times C_2$ that identifies the two copies of $u$ and the two copies of $v$ is injective on $S$, so by Lemma~\ref{lem:quot} we get a counterexample in
\[ \G=A_6\times A_6\times\langle u\rangle\times\langle v\rangle,\qquad|\G|=518\,400. \]
Explicitly, with $\epsilon(h)\in\{0,1\}$ equal to $1$ exactly on the three involutions of $H_i$, and writing $\bar h$ for the $A_6$-part of $h$:
\begin{align*}
A&=\{(\bar h,\bar h',u^{\epsilon(h)+\epsilon(h')},1):h,h'\in H_1\}, & D&=\{(\bar g,\bar g',1,v^{\epsilon(g)+\epsilon(g')}):g,g'\in H_2\},\\
B&=\{(\bar g,e,1,v^{\epsilon(g)}),\,(e,\bar g,1,v^{\epsilon(g)}):g\in H_2^\#\}, & C&=\{(\bar h,e,u^{\epsilon(h)},1),\,(e,\bar h,u^{\epsilon(h)},1):h\in H_1^\#\}.
\end{align*}
These have $|A|=|D|=36$ and $|B|=|C|=10$. Direct enumeration gives
\[ |S|=66096,\quad N_1=N_4=2160,\quad N_2=N_3=7776, \]\[ \frac{|S|^3}{N_1N_2N_3N_4}=\frac{132651}{129600}=\frac{51^3}{4\cdot25\cdot36^2}\approx1.0235. \]
Both central factors are needed. Identifying $u$ with $v$, or dropping either of them, destroys the example: the ratio drops to $0.27$ or $0.39$. Without the twists, $A_6\times A_6$ has no configuration of this type at all, because $|A_6|=360<936=kl(1+rs)$.
\end{example}

In lowest terms the ratio is $4913/4800$. The same numbers arise in $S_6\times S_6\times C_2$ (order $1\,036\,800$) from $H_1=\langle(1\,2\,3\,4\,5\,6)\rangle$ and a twisted $S_3$.

\subsection{Alternating groups}

\begin{example}[$A_6^3$: two $5$-cycles]
In $A_6$ let $\sigma=(2\,3\,4\,5\,6)$ and $\tau=(1\,2\,3\,4\,5)$, and put $H_1=\langle\sigma\rangle$, $H_2=\langle\tau\rangle$, $Z'=\{\sigma,\sigma^2,\sigma^4\}$, $Y'=\{\tau,\tau^2,\tau^3\}$. The $10$ double cosets $H_1H_2$ and $H_1\tau^i\sigma^jH_2$ are distinct and free. With $m=3$ coordinates the sets $A=H_1^3$, $D=H_2^3$, $B=\{y^{(i)}\}$, $C=\{z^{(j)}\}$ in $A_6^3$ (sizes $125,9,9,125$) give
\[ |S|=437500,\quad(N_1,N_2,N_3,N_4)=(4500,62500,62500,4500),\qquad\frac{|S|^3}{N_1N_2N_3N_4}\approx1.0586, \]
and with $m=4$ in $A_6^4$ the ratio is $\approx1.3741$. Both values were checked by brute force. This is perhaps the simplest counterexample of all: two $5$-cycles in $A_6$.
\end{example}
\begin{example}[$A_5^4$]\label{ex:a54}
In $G_0=A_5\times A_5$ (order $3600$) take two ``diagonal'' subgroups:
\begin{align*}
H_1&=\langle(1\,2\,4\,5\,3)(6\,8\,10\,9\,7)\rangle\cong C_5,\\
H_2&=\langle(1\,5\,4\,2\,3)(6\,10\,8\,9\,7),\ (2\,4)(3\,5)(6\,10)(7\,8)\rangle\cong D_{10},
\end{align*}
where the second $A_5$ acts on $\{6,\dots,10\}$. Then $\langle H_1,H_2\rangle=G_0$. With $Z'=H_1^\#$ ($s=4$) and $Y'=H_2^\#\setminus\{(1\,2)(4\,5)(6\,7)(9\,10)\}$ ($r=8$), the $33$ double cosets $H_1H_2$ and $H_1yzH_2$ are distinct and free. Theorem~\ref{thm:multi} with $m=2$ then gives a counterexample in $G_0\times G_0=A_5^4$ (order $12\,960\,000$). Its sets have sizes $|A|=25$, $|B|=16$, $|C|=8$, $|D|=100$. Direct enumeration gives
\[ |S|=162500,\quad(N_1,N_2,N_3,N_4)=(7200,12500,22500,2000), \]\[ \frac{|S|^3}{N_1N_2N_3N_4}=\frac{65^3}{4\cdot32\cdot81\cdot25}\approx1.0595. \]
\end{example}
\begin{example}[$A_8\times A_8$: a large certified ratio]\label{ex:a8}
In $A_8$ let $H_1\cong C_4\times C_2$ be the regular subgroup
\[ H_1=\langle(1\,2\,8\,7)(3\,6\,5\,4),\ (1\,3)(2\,6)(4\,7)(5\,8)\rangle, \]
and let $H_2\cong C_3\rtimes C_4$ be the point stabiliser subgroup
\[ H_2=\langle(6\,7\,8),\ (2\,3\,4\,5)(7\,8)\rangle. \]
Explicitly, $H_1=\{e$, $(1\,2\,8\,7)(3\,6\,5\,4)$, $(1\,8)(2\,7)(3\,5)(4\,6)$, $(1\,7\,8\,2)(3\,4\,5\,6)$, $(1\,3)(2\,6)(4\,7)(5\,8)$, $(1\,4\,8\,6)(2\,3\,7\,5)$, $(1\,5)(2\,4)(3\,8)(6\,7)$, $(1\,6\,8\,4)(2\,5\,7\,3)\}$. These are the eight elements of $H_1$, and the displayed generators of $H_1$ and $H_2$ generate $A_8$. Take $Z'=H_1^\#$ ($s=7$) and $Y'=H_2^\#\setminus\{(2\,5\,4\,3)(6\,7)\}$ ($r=10$). The $71$ double cosets $H_1H_2$ and $H_1yzH_2$ are distinct and free. In $A_8\times A_8$ the sets have sizes $64,20,14,144$, and direct enumeration gives
\[ |S|=1\,299\,456,\quad(N_1,N_2,N_3,N_4)=(22176,\,73728,\,101376,\,10240), \]
\[ \frac{|S|^3}{N_1N_2N_3N_4}=\frac{141^3}{4\cdot70\cdot121\cdot64}\approx1.2928. \]
This is a comparatively strong, explicitly displayed and independently certified four-set example. It is not the largest ratio reported here: Table~\ref{tab:alt} lists ratios of about $1.430$ ($A_{10}^2$) and $1.555$ ($A_{12}^2$); the latter rests on sampled subgroups only. The search of $A_8$ that found this pair reported $28\,831$ rectangles; that global enumeration was not independently rerun.
This is by far the largest factor we have found. A complete search of $A_8$ (subgroups of order $5$--$12$, double cosets computed directly, without DoubleCosets) finds $28831$ good rectangles. All pairs generate $A_7$ or $A_8$; in $A_7$ the best is $(r,s)=(7,4)$ with ratio $1.0334$, giving a counterexample in $A_7\times A_7$ (order $6.35\cdot10^6$). The groups are large ($|A_8|^2\approx4.1\cdot10^8$), so these examples do not compete with $A_6^2\times C_2^2$ in size.
\end{example}
\subsection{Larger factors}\label{sec:larger}
For fixed per-coordinate data $(r,s)$ the factor grows linearly in the number $m$ of coordinates, like $m\bigl(rs/((r+1)(s+1))\bigr)^2$. The best rectangles found in alternating groups are listed in Table~\ref{tab:alt}.
\begin{table}[htbp]\begin{center}
\begin{tabular}{lllllll}
\hline
$G_0$ & $(H_1,H_2)$ & $(r,s)$ & $m=2$ & $m=3$ & $m=4$ & $m=6$\\
\hline
$A_6$ & $C_5,C_5$ & $(3,3)$ & $0.744$ & $\mathbf{1.059}$ & $\mathbf{1.374}$ & $2.006$\\
$A_7$ & $C_5,D_{12}$ & $(7,4)$ & $1.033$ & $1.523$ & $2.013$ & $2.993$\\
$A_8$ & $C_4\times C_2,\ C_3\rtimes C_4$ & $(10,7)$ & $\mathbf{1.293}$ & $1.926$ & $2.558$ & $3.824$\\
$A_9$ & $C_6\times C_2,\ D_{12}$ & $(9,10)$ & $1.361$ & $2.031$ & $2.700$ & $4.039$\\
$A_{10}$ & $C_6\times C_2,\ C_6\times C_2$ & $(11,11)$ & $\mathbf{1.430}$ & $2.136$ & $2.842$ & $4.254$\\
$A_{12}$ & two subgroups of order $16$ & $(15,15)$ & $1.555$ & $2.328$ & $3.100$ & $4.645$\\
\hline
\end{tabular}
\end{center}
\caption{Best rectangles $(r,s)$ found in alternating groups, and the resulting ratios of Theorem~\ref{thm:multi} for $m=2,3,4,6$ coordinates. Bold entries were confirmed by enumeration of the actual sets.}\label{tab:alt}\end{table}
The $A_6$, $A_7$ and $A_8$ rows have local double-coset certificates, and the $A_6^3$, $A_6^4$, $A_8^2$ and $A_{10}^2$ ratios were confirmed by enumeration of the actual sets. The $A_9$, $A_{10}$ and $A_{12}$ pairs were found by random sampling of subgroups.

\subsection{Further examples}

\begin{example}[$S_7\times S_7$, cyclic subgroups]\label{ex:s7}
Let $\sigma=(1\,2\,3\,4\,5\,6\,7)$, $\tau=(1\,6\,7\,2\,3\,4)$, $H_1=\langle\sigma\rangle\cong C_7$ and $H_2=\langle\tau\rangle\cong C_6$. The $31$ double cosets $H_1H_2$ and $H_1\tau^i\sigma^jH_2$ are distinct and free ($r=5$, $s=6$). With $A=H_1^2$, $D=H_2^2$, $B=(H_2^\#\times1)\cup(1\times H_2^\#)$ and $C=(H_1^\#\times1)\cup(1\times H_1^\#)$ one gets $|S|=107604$, $(N_1,N_2,N_3,N_4)=(2592,12348,10584,3430)$ and ratio $226981/211680\approx1.0723$. This is the largest ratio among our small examples, in a group of order $2.5\cdot10^7$.
\end{example}
Twisting by non-central factors also works: in $(A_5\times S_3\times D_{10})^2$ a twisted $S_3$ and a twisted $D_{10}$ give $(r,s)=(6,4)$ and ratio $117649/117600$, and $A_5\times D_{10}\times D_{10}$ gives $(4,9)$ and ratio $1.0806$; none of these groups is smaller than $A_6^2\times C_2^2$.

\begin{table}[htbp]\begin{center}
\begin{tabular}{lllr}
\hline
group $\G$ & $H_1,\ H_2$ & $(r,s)$, $m$ & $|\G|$\\
\hline
$A_6\times A_6\times C_2^2$ & twisted $S_3$, twisted $S_3$ & $(5,5)$, $2$ & $518\,400$\\
$S_6\times S_6\times C_2$ & $C_6$, twisted $S_3$ & $(5,5)$, $2$ & $1\,036\,800$\\
$\mathrm{PGL}(2,13)^2$ & $C_6$, $S_3$ & $(5,5)$, $2$ & $4\,769\,856$\\
$\mathrm{PSL}(2,17)^2$ & $S_3$, $S_3$ & $(5,5)$, $2$ & $5\,992\,704$\\
$(A_5\times S_3\times D_{10})^2$ & twisted $S_3$, twisted $D_{10}$ & $(6,4)$, $2$ & $12\,960\,000$\\
$A_5^4$ & diagonal $C_5$, diagonal $D_{10}$ & $(8,4)$, $2$ & $12\,960\,000$\\
$A_7^2$ & $C_5$, $D_{12}$ & $(7,4)$, $2$ & $6\,350\,400$\\
$A_8^2$ (ratio $1.29$) & $C_4\times C_2$, $C_3\rtimes C_4$ & $(10,7)$, $2$ & $406\,425\,600$\\
$S_7^2$ & $C_7$, $C_6$ & $(5,6)$, $2$ & $25\,401\,600$\\
$A_6^3$ & $C_5$, $C_5$ & $(3,3)$, $3$ & $46\,656\,000$\\
class-two $2$-group (Thm.~\ref{thm:c2}) & $V$, $V$ & sparse $9\times9$ pattern, $1$ & $2^{23}\approx8.4\cdot10^{6}$\\
$D_4^{13}$ & $C_2$, $C_2$ & $(1,1)$, $13$ & $\approx5.5\cdot10^{11}$\\
$D_4^{14}$ [NNZ] & & & $\approx4.4\cdot10^{12}$\\
\hline
\end{tabular}
\end{center}
\caption{Four-set counterexample groups. Ratios $>1$ are certified by Theorem~\ref{thm:multi} (or Theorem~\ref{thm:c2}); actual-set enumerations cover $A_6^2\times C_2^2$, $S_6^2\times C_2$, $A_5^4$, $A_8^2$, $S_7^2$, $A_6^3$ and the class-two group. Rows based on reported searches are qualified in Section~\ref{sec:method}.}\label{tab:witness}\end{table}
\subsection{Torsion-free groups}\label{sec:tf}
Theorem~\ref{thm:all2} gives counterexamples in direct powers of torsion-free non-abelian groups. A two-coordinate example already exists in $F_2\times F_2$, where $F_2=\langle a,b\rangle$ is the free group on two generators. It is the construction of Theorem~\ref{thm:multi} with the finite subgroups replaced by segments of infinite cyclic groups. Fix integers $k\ge2$ and $r,s\ge1$, write $P=\{b^u:0\le u<k\}$ and $Q=\{a^v:0\le v<k\}$, and in $F_2\times F_2$ take
\[ A=P\times P,\quad D=Q\times Q,\quad B=\{(a^i,e),(e,a^i):1\le i\le r\},\quad C=\{(b^j,e),(e,b^j):1\le j\le s\}. \]
\begin{proposition}\label{prop:f2}
With these sets,
\begin{align*}
|ABCD|&=2k^2(k+r-1)(k+s-1)-(k-1)^4+2rsk^4,\\
N_1&=2sk(k+r-1)+2rsk^2, & N_4&=2rk(k+s-1)+2rsk^2,\\
N_2&=sk^4+k^3(k+s-1), & N_3&=rk^4+k^3(k+r-1).
\end{align*}
\end{proposition}
The formulas come from reducing words in the free group: the first term of $|ABCD|$ counts the two ``cross'' orientations of reduced-word rectangles, $(k-1)^4$ is their overlap, and the remaining $2rsk^4$ elements have two distinct mixed-syllable orientations. They were checked against exhaustive enumeration for $k\le5$. For $k=150$ and $r=s=5$,
\[ |ABCD|=25\,886\,835\,599,\qquad (N_1,N_2,N_3,N_4)=(1\,356\,000,\ 3\,051\,000\,000,\ 3\,051\,000\,000,\ 1\,356\,000), \]
so $R\approx1.01352>1$. As $k\to\infty$ the ratio tends to $(1+2rs)^3/(4rs(r+1)^2(s+1)^2)$, which is $51^3/(4\cdot25\cdot36^2)\approx1.0235$ for $r=s=5$. Hence the inequality fails in every group containing $F_2\times F_2$, for instance $\mathrm{SL}_4(\mathbb Z)$ and $\mathrm{Aut}(F_n)$ for $n\ge4$.

\subsection{A class-two $2$-group}\label{sec:c2}
The next four-set example lives in a single nilpotent group rather than in a direct power. It is another instance of the double-coset collapse, now with a sparse one-coordinate incidence pattern.

Let $V=\mathbb F_2^9$ with basis $e_0,\dots,e_8$ (also written $f_0,\dots,f_8$ for the second factor), $W=\mathbb F_2^5$, and let $\omega:V\times V\to W$ be bilinear. Among the $9\times9$ values $\omega(e_i,f_j)$, assign distinct nonzero vectors of $W$ to the $27$ positions with $j-i\equiv0,1,2\pmod9$, and zero to the remaining $54$ positions. Define a multiplication on $G=V\times V\times W$ by
\[ (c,b,w)\,(c',b',w')=(c+c',\ b+b',\ w+w'+\omega(b,c')). \]
Then $|G|=2^{23}$ and $G$ is nilpotent of class two. Put
\begin{align*} H&=\{(c,0,0):c\in V\}, & J&=\{(0,b,0):b\in V\},\\ B&=\{e\}\cup\{(0,e_i,0):0\le i<9\}, & C&=\{e\}\cup\{(f_j,0,0):0\le j<9\}. \end{align*}

\begin{theorem}\label{thm:c2}
For the four sets $(H,B,C,J)$,
\[ |HBCJ|=28\cdot2^{18}=7\,340\,032,\qquad (N_1,N_2,N_3,N_4)=(37\cdot2^9,\ 2^{20},\ 2^{20},\ 37\cdot2^9), \]
\[ \frac{|HBCJ|^3}{N_1N_2N_3N_4}=\frac{28^3}{37^2\cdot16}=\frac{1372}{1369}>1. \]
\end{theorem}
\begin{proof}
$H$ and $J$ are subgroups of order $2^9$ with $H\cap J=\{e\}$, and $HJ=\{(c,b,0)\}$. For $y=(0,e_i,0)\in B$ and $z=(f_j,0,0)\in C$ we have $yz=zy\,(0,0,\omega(e_i,f_j))$, and $HyzJ=V\times V\times\{\omega(e_i,f_j)\}$ is a full double coset of size $2^{18}$. The $27$ incident pairs give the $27$ distinct nonzero labels, and all other pairs give $HJ$. So $S=HBCJ$ is the disjoint union of $28$ double cosets of size $2^{18}$. Every row and every column of the incidence pattern contains three incidences, so fixing $y$ or $z$ leaves $1+3=4$ double cosets: $N_2=N_3=4\cdot2^{18}$. In $BCJ$, the right $J$-cosets are indexed by the identity, the nine elements of $C^\#$ (giving $f_jJ$) and the $27$ incident pairs (giving $(f_j,0,\omega(e_i,f_j))J$); non-incident pairs are absorbed into $f_jJ$. Hence $N_1=37\cdot2^9$, and symmetrically $N_4=37\cdot2^9$.
\end{proof}

This example also verifies the following general template, which does not require a full rectangle and works with one coordinate.

\begin{proposition}[Sparse single-coordinate template]\label{prop:sparse}
Let $H,J\le G$ with $|H|=k$, $|J|=\ell$, $H\cap J=\{e\}$, let $Y\subseteq J^\#$ and $Z\subseteq H^\#$ have sizes $r,s$, and let $E\subseteq Y\times Z$ with $|E|=I$. Assume that (i) for $(y,z)\in E$ the $I$ double cosets $HyzJ$ are full, pairwise distinct and distinct from $HJ$, and (ii) for $(y,z)\notin E$ the elements $y,z$ commute. For $A=H$, $D=J$, $B=\{e\}\cup Y$, $C=\{e\}\cup Z$, and with $\Delta_Y,\Delta_Z$ the maximal row and column degrees of $E$,
\[ |ABCD|=k\ell(1+I),\quad N_2=k\ell(1+\Delta_Y),\quad N_3=k\ell(1+\Delta_Z),\quad N_1=\ell(1+s+I),\quad N_4=k(1+r+I), \]
\[ R=\frac{(1+I)^3}{(1+s+I)(1+r+I)(1+\Delta_Y)(1+\Delta_Z)} . \]
\end{proposition}
\begin{proof}
The full double cosets count $S$, and the row and column degrees count the two middle slices. In $BCJ$, a non-incident product $yz=zy$ is absorbed into $zJ$, so the right $J$-cosets are $J$, the $s$ cosets $zJ$, and the $I$ cosets $yzJ$ for $(y,z)\in E$; the latter are distinct and new by (i). This gives $N_1$; $N_4$ is symmetric.
\end{proof}

Hypothesis (ii) cannot be weakened to $yz\in HJ$: in $S_3=\langle a\rangle\rtimes\langle s\rangle$ with $Y=\{s\}$, $Z=\{a\}$ one has $sa\in HJ$ but $saJ=a^{-1}J\ne aJ$. For a $d$-regular pattern on $t+t$ vertices with all $td$ labels distinct the ratio is $(1+td)^3/((1+t+td)^2(1+d)^2)$: $(t,d)=(9,3)$ is Theorem~\ref{thm:c2}, and $(t,d)=(12,3)$ gives $37^3/(49^2\cdot16)\approx1.318$ in a class-two group of order $2^{30}$. This sparse template lies outside the full-rectangle bounds of Lemma~\ref{lem:derived} and Corollary~\ref{cor:size}; whether it can give a four-set witness in a group smaller than $A_6^2\times C_2^2$ is open.

\section{Other values of $n$}\label{sec:n}

\subsection{Small $n$}
For $n=2$ the inequality reads $|A_1A_2|\le N_1N_2=|A_2|\,|A_1|$, which is trivial. For $n=3$ it holds in every group: Ruzsa's proof of $|A+B+C|^2\le|A+B||B+C||A+C|$ carries over to groups in the form $|ABC|^2\le|BC|\cdot\max_{b\in B}|AbC|\cdot|AB|$ ([GMR], remark before Problem~1.3; [NNZ]). So $n=4$ is the first value where the inequality can fail, and it does.

\subsection{All $n\ge4$, including non-degenerate examples}
Lemma~\ref{lem:pad} turns every counterexample for $n=4$ into one for every $n\ge4$ by adding singletons. Singletons can be avoided. Let $(A_1,\dots,A_4)$ be a counterexample in $\G$. In $\G\times C_2^{\,n-4}$ take $A_i\times\{0\}$ for $i\le4$ and $\{e\}\times\{0,e_j\}$ for $5\le j\le n$. The abelian part is a ``box'', for which \eqref{eq:gmr} is an equality: it multiplies $|S|^{n-1}$ and $\prod N_i$ by the same power $2^{(n-4)(n-1)}$. The ratio is therefore unchanged, and all $n$ sets have at least two elements. Consequently the constructions of Theorems~\ref{thm:binary}--\ref{thm:multi}, Corollary~\ref{cor:iff}, and the $S_7\times S_7$ and $F_2\times F_2$ examples extend to every $n\ge4$. In particular, \emph{for every $n\ge4$, the inequality \eqref{eq:nagmr} holds in all direct powers of $G$ if and only if $G$ is abelian}.

\subsection{How badly can it fail? The best exponent}\label{sec:exp}
Ratios are multiplicative under direct products. If $(A_i)$ in $\G$ and $(A_i')$ in $\G'$ are given, then in $\G\times\G'$ the sets $A_i\times A_i'$ have $|S|=|S||S'|$ and $N_i=N_iN_i'$, because a slice of a product set is the product of the slices. Hence one counterexample with ratio $R>1$ gives ratio $R^t$ in the $t$-th power. So for $n\ge4$ \emph{no constant} $C$ makes $|S|^{n-1}\le C\prod N_i$ true. The right question concerns the exponent. Let
\[ \theta_n=\inf\ \frac{\log\prod_iN_i}{\log|S|}, \]
the infimum over all groups and all finite $A_1,\dots,A_n$ with $|S|\ge2$. Then $\theta_n=n-1$ would be the conjecture (with constant $1$).

\begin{proposition}\label{prop:exp}
$\theta_3=2$, and for $n\ge4$
\[ 2\ \le\ \theta_n\ \le\ n-1-0.0768. \]
More precisely $\theta_4\le2.9232$, and $\theta_n\le\theta_4+(n-4)$.
\end{proposition}

\begin{proof}
\emph{Lower bound.} Apply the case $n=3$ to the three blocks $A_1\cdots A_{i-1}$, $A_i$, $A_{i+1}\cdots A_n$ ($2\le i\le n-1$). This gives $|S|^2\le|A_i\cdots A_n|\,N_i\,|A_1\cdots A_i|\le N_1N_iN_n$, since $A_i\cdots A_n$ is contained in a translate of $A_2\cdots A_n$, and similarly on the other side. As all $N_j\ge1$, $|S|^2\le\prod_jN_j$.

\emph{Upper bound for $n=4$.} By the type-class argument of [NNZ, Prop.~B.1], an entropic configuration gives sets in $G^m$ with $\frac1m\log|S_m|\to h:=\max H(U)$ and $\frac1m\log N_{i,m}\to C_i:=\max H(U\mid X_i)$. Hence $\theta_4\le\sum_iC_i/h$. For the multi-rare configuration of Theorem~\ref{thm:multi} with $X_1,X_4$ uniform on $H_1,H_2$, $Y\in\{e\}\cup Y'$, $Z\in\{e\}\cup Z'$, $\Pr(\text{rare})=\varepsilon$ and the rare values uniform, a conditional-entropy count of the disjoint double-coset outcomes gives
\begin{align*}
h&=\log kl+h_2(\varepsilon)+\varepsilon\log rs,\\
\textstyle\sum_iC_i&\le 3\log kl+2h_2(\varepsilon)+\varepsilon\bigl(\log rs+2\log(r+1)(s+1)\bigr).
\end{align*}
For the configuration of Example~\ref{ex:a6} ($k=l=6$, $r=s=5$) and $\varepsilon=0.258$ the quotient is $2.9232$ (for the $D_4$ configuration of [NNZ] it is $2.9622$). The bound for $n\ge5$ follows by adding singletons, which adds $(n-4)\log|S|$ to $\log\prod N_i$.
\end{proof}

So the deficit $(n-1)-\theta_n$ is positive for all $n\ge4$. The four-set constructions do not make it grow with $n$; the family of Section~\ref{sec:fam} does, up to about $0.356$ (Corollary~\ref{cor:theta}).

\subsection{The case $n=6$: a five-double-coset criterion}\label{sec:a5n6}
Six factors permit much smaller groups. The criterion below is local: it concerns only five double cosets in each coordinate.

\begin{theorem}[Five double cosets]\label{thm:crit}
Let $G_1,G_2$ be groups and $a_i,b_i\in G_i$ elements of the same finite order $p_i\ge3$; 
put $H_i=\langle a_i\rangle$, $K_i=\langle b_i\rangle$. Suppose that in each $G_i$ the five sets
\begin{equation}\label{eq:five}
H_iK_i,\qquad H_ib_i^{\epsilon}a_i^{\delta}K_i\quad(\epsilon,\delta\in\{\pm1\})
\end{equation}
are pairwise distinct double cosets of size $p_i^2$. In $G_1\times G_2$ let
\begin{align*}
A_1&=H_1\times H_2, & A_2&=\{(b_1^{-1},e),(e,b_2)\}, & A_3&=\{(b_1,e),(e,b_2^{-1})\},\\
A_4&=\{e,(a_1^{-1},e),(e,a_2^{-1})\}, & A_5&=\{e,(a_1,e),(e,a_2)\}, & A_6&=K_1\times K_2.
\end{align*}
With $P=p_1p_2$, $|S|=13P^2$ and $(N_1,\dots,N_6)=(19P,7P^2,7P^2,5P^2,5P^2,15P)$, so
\[ \frac{|S|^5}{N_1\cdots N_6}=\frac{13^5}{19\cdot7^2\cdot5^2\cdot15}=\frac{371293}{349125}\approx1.0635 . \]
\end{theorem}

\begin{proof}
$A_2A_3=\{e,(b_1,b_2),(b_1^{-1},b_2^{-1})\}$, and $A_4A_5=\{(a_1^s,a_2^t):(s,t)\in T\}$ with $T=\{(0,0),(\pm1,0),(0,\pm1),(1,-1),(-1,1)\}$. Hence $S$ is the union of the boxes $H_1b_1^ua_1^sK_1\times H_2b_2^ua_2^tK_2$ ($u\in\{0,\pm1\}$, $(s,t)\in T$). For $u=0$, and for $(s,t)=(0,0)$, this is the base box $H_1K_1\times H_2K_2$. For $u=\pm1$ the other six pairs $(s,t)$ give boxes with a non-trivial double coset in at least one coordinate. By \eqref{eq:five} these $1+6+6=13$ boxes are pairwise disjoint, and each has $P^2$ elements.

Fixing an element of $A_2$ or $A_3$ leaves one value of $u\ne0$ besides $u=0$, hence $1+6=7$ boxes. Fixing an element of $A_4$ or $A_5$ leaves three pairs $(s,t)$ (for instance $(0,0),(1,0),(0,1)$), hence $1+2+2=5$ boxes. For $N_1=|A_2\cdots A_6|$ we count right cosets $b_1^ua_1^sK_1\times b_2^ua_2^tK_2$ of size $P$. For $u=0$ the seven pairs in $T$ give seven distinct cosets, since $H_i\cap K_i=\{e\}$. For $u=\pm1$, $b_i^ua_i^sK_i=K_i$ if $s=0$, and otherwise it lies in a non-trivial double coset; this gives $6+6$ new cosets. So $N_1=19P$. For $N_6=|A_1\cdots A_5|$ we count left cosets $H_1b_1^ua_1^s\times H_2b_2^ua_2^t$: for $u=0$ all seven collapse to $H_1\times H_2$, while for $u=\pm1$ the seven pairs are distinct, so $N_6=(1+7+7)P$. All other choices give the same counts by translation.
\end{proof}

\begin{corollary}\label{thm:a5n6}
For $a=(3\,4\,5)$ and $b=(1\,2\,3)$ in $A_5$, the sets of Theorem~\ref{thm:crit} in $A_5\times A_5$ (order $3600$) have $|S|=1053$, $(N_i)=(171,567,567,405,405,135)$ and ratio $371293/349125$.
\end{corollary}
\begin{proof}
$K=\langle b\rangle$ fixes $4$ and $5$, so the right coset $gK$ is determined by $(g(4),g(5))$. Left multiplication by $H=\langle a\rangle$ rotates both entries by $(3\,4\,5)$. The five representatives $e,\ b^{\pm1}a^{\pm1}$ give five different $H$-orbits of ordered pairs, each of size $3$, so \eqref{eq:five} holds with $p=3$.
\end{proof}

\begin{theorem}[Characteristic $3$]\label{thm:char3}
Let $q=3^f$ with $f\ge2$, $t\in\mathbb F_q\setminus\mathbb F_3$, and let $a,b$ be the images of $\begin{pmatrix}1&1\\0&1\end{pmatrix}$ and $\begin{pmatrix}1&0\\t&1\end{pmatrix}$ in $\mathrm{PSL}(2,q)$. Then $a,b$ satisfy \eqref{eq:five} with $p=3$. Hence $\mathrm{PSL}(2,q)\times G_2$ has a six-set counterexample whenever $G_2$ satisfies \eqref{eq:five} (for instance $G_2=A_5$ or $\mathrm{PSL}(2,q')$, $q'=3^{f'}$).
\end{theorem}
\begin{proof}
$H=\langle a\rangle$ consists of upper and $K=\langle b\rangle$ of lower unitriangular matrices, so $H\cap K=\{e\}$. Left multiplication by $H$ only changes the first row, and right multiplication by $K$ only changes the first column. So the bottom-right entry $d$ of a representative is an invariant of the double coset, up to sign. For $x=b^ia^j=\begin{pmatrix}1&j\\ it&1+ijt\end{pmatrix}$ we get $d=1$ ($x=e$), $1+t$ ($ij=1$), $1-t$ ($ij=-1$). These are pairwise different up to sign, since $t\notin\mathbb F_3$ and $2\ne0$.

To separate $(i,j)$ from $(-i,-j)$, note that the bottom rows of the elements of $HxK$ are $\pm(it+cdt,\ d)$ with $c\in\mathbb F_3$. Suppose $x=b^ia^j$ and $x'=b^{-i}a^{-j}$, which have the same $d\ne0$, lie in the same double coset. The sign must be $+$, since $d=-d$ would force $d=0$. Then $it+cdt=-it+c'dt$, i.e.\ $2it=(c'-c)dt$; if $c=c'$ this gives $2it=0$, impossible, and otherwise $d=2i/(c'-c)\in\mathbb F_3$, contradicting $d=1\pm t\notin\mathbb F_3$.

Finally, $|HxK|=9$ exactly when $H\cap xKx^{-1}=\{e\}$. The non-trivial elements of $H$ fix only the point $\infty$ of the projective line, and those of $xKx^{-1}$ fix only $x(0)$. Since $d\ne0$ we have $x(0)\ne\infty$, so the intersection is trivial.
\end{proof}

The criterion has been checked (by computer) for $\mathrm{PSL}(2,7)$, with $a=(2\,4\,6)(3\,5\,7)$ and $b=(1\,2\,3)(5\,6\,7)$ on the points of the Fano plane, for $\mathrm{PSL}(2,8)$, and for $C_3\times\mathrm{SL}(2,3)$ (Section~\ref{sec:1728}). For the Suzuki group $\mathrm{Sz}(8)$ an explicit pair of elements of order $5$ satisfies it; the matrices and the verification are given in Appendix~\ref{app:sz}.

\begin{remark}[What condition \eqref{eq:five} does and does not give]\label{rem:five}
Condition \eqref{eq:five} in $G_1$ and in $G_2$ yields a six-set counterexample in $G_1\times G_2$, hence by Proposition~\ref{prop:sub} in every group containing $G_1\times G_2$, and by Lemma~\ref{lem:pad} for every $n\ge6$. It says nothing about $G_1$ alone: $A_5$ satisfies \eqref{eq:five}, but no counterexample in $A_5$ itself is known. Nor does its failure say anything positive: $D_4$ fails \eqref{eq:five}, yet $D_4^{13}$ is a counterexample for $n=4$ by Theorem~\ref{thm:binary}. The condition is restrictive for small groups. It requires two cyclic subgroups of the same order $p\ge3$ meeting trivially, so dihedral and dicyclic groups and all groups with a unique subgroup of each order fail it; and the five disjoint double cosets force $|\langle a,b\rangle|\ge5p^2\ge45$ and $|\langle a,b\rangle'|\ge5$. An exhaustive computer check over all pairs $a,b$ shows that among the $1017$ non-abelian groups of order at most $120$ exactly seven satisfy \eqref{eq:five}: $A_5$, $C_3\times\mathrm{SL}(2,3)$, $C_3\wr C_3$, \texttt{SmallGroup(96,3)}, $\mathrm{SL}(2,5)$, $S_5$ and $C_2\times A_5$, all with $p=3$.

For simple groups the picture is different: every non-abelian finite simple group we tested satisfies \eqref{eq:five} (Question~\ref{q:simple} lists them). Consequently all products $G_1\times G_2$ of two of these simple groups have six-set counterexamples. Some simple groups contain counterexamples themselves: the ten-set configuration of Theorem~\ref{thm:243} is realised in $\mathrm{PSL}(2,8)$, $\mathrm{PSL}(2,11)$ and $\mathrm{PSL}(2,13)$, the fourteen-set member of the family of Theorem~\ref{thm:fam} in $\mathrm{PSL}(2,17)$, $\mathrm{PSL}(2,19)$ and (with ratio $1.12$) in $A_7$, and Corollary~\ref{cor:famfinite} gives counterexamples in $\mathrm{PSL}(2,p)$ for every prime $p>2^{2s+3}$. By subgroup inheritance, $A_{10}\supseteq A_5\times A_5$ has a six-set counterexample, and $A_{14}\supseteq A_7\times A_7$, $A_{16}\supseteq A_8\times A_8$ and $A_{18}\supseteq A_6^3$ have four-set counterexamples.
\end{remark}

\emph{How it was found.} We extended the counting of Theorem~\ref{thm:multi} to a ``generic'' model. In it, $\langle a\rangle\cong\mathbb Z_k$ and $\langle b\rangle\cong\mathbb Z_l$ generate their free product, $A_1=H_1^m$ and $A_n=H_2^m$, and the middle sets consist of single letters in single coordinates. All cardinalities are then counts of reduced words modulo $H_1$ and $H_2$. Simulated annealing over the middle sets produced the configuration above ($n=6$, $m=2$, $k=l=3$) and the one of Theorem~\ref{thm:243} ($n=10$, $m=1$).

\subsection{A six-set counterexample of order $1728$}\label{sec:1728}
Theorem~\ref{thm:crit} only needs the two coordinates to commute. By Lemmas~\ref{lem:quot} and~\ref{lem:twist} it therefore applies in central products whenever the quotient map is injective on $S$.

\begin{theorem}\label{thm:1728}
In $\mathrm{GL}(2,7)$ let
\[ a=\begin{pmatrix}5&5\\5&1\end{pmatrix},\qquad b=\begin{pmatrix}6&5\\5&6\end{pmatrix},\qquad w=\begin{pmatrix}2&0\\0&2\end{pmatrix}. \]
Then $G_0=\langle a,b\rangle\cong C_3\times\mathrm{SL}(2,3)$ has order $72$, the element $w\in Z(G_0)$ has order $3$, and $a,b$ satisfy \eqref{eq:five}. In $G=(G_0\times G_0)/\langle(w,w)\rangle$, a group of order $1728$, the sets of Theorem~\ref{thm:crit} have $|S|=1053$ and ratio $13^5/349125$.
\end{theorem}
\begin{proof}
Condition \eqref{eq:five} is a finite check in $G_0$, so Theorem~\ref{thm:crit} applies in $G_0\times G_0$. Represent each coset of $\langle(w,w)\rangle$ by the lexicographically least member of its three-element orbit; the $1053$ elements of $S$ have pairwise distinct representatives, so the quotient map is injective on $S$, and Lemma~\ref{lem:quot} applies.
\end{proof}

The anti-diagonal $\langle(w,w^{-1})\rangle$ works too. A quotient of order $864<1053$ is impossible. Among all groups $G_0$ of order at most $503$ generated by elements of order $3$, the smallest for which the configuration works in $G_0\times G_0$ are $A_5$, $C_3\times\mathrm{SL}(2,3)$, $C_3\wr C_3$ (central quotient of order $2187$), and groups of orders $96$--$180$. No mixed product $A_5\times X$ or $(C_3\times\mathrm{SL}(2,3))\times X$ with $|X|<60$ works.

\subsection{Seven and eight sets without singletons}\label{sec:n8}
Padding produces singleton sets. Genuine $7$- and $8$-set examples with \emph{larger} ratio come from factorising the outer subgroups. For $g$ of order $3$ put
\[ L_+(g)=\{e,(g,g),(g^2,g^2)\},\qquad L_-(g)=\{e,(g,g^2),(g^2,g)\},\qquad L_-(g)L_+(g)=\langle g\rangle\times\langle g\rangle. \]

\begin{theorem}\label{thm:n8}
Let $a,b$ be as in Corollary~\ref{thm:a5n6} (in $A_5\times A_5$) or Theorem~\ref{thm:1728} (in the group of order $1728$), and let $A_2,\dots,A_5$ be the middle sets of Theorem~\ref{thm:crit}.
\begin{enumerate}
\item The $8$ sets $L_-(a),L_+(a),A_2,A_3,A_4,A_5,L_-(b),L_+(b)$ have $|S|=1053$ and $(N_i)=(405,405,567,567,405,405,351,405)$.
\item The $7$ sets $L_-(a),L_+(a),A_2,A_3,A_4,A_5,\langle b\rangle\times\langle b\rangle$ have $|S|=1053$ and $(N_i)=(405,405,567,567,405,405,135)$.
\end{enumerate}
In both cases the ratio is
\[ \frac{1053^7}{405^5\cdot351\cdot567^2}=\frac{1053^6}{405^4\cdot567^2\cdot135}=\frac{4826809}{4134375}\approx1.1675 , \]
and no set is a singleton.
\end{theorem}
\begin{proof}
The outer sets multiply to $\langle a\rangle^2$ and $\langle b\rangle^2$, so $S$ and the middle slices are those of Theorem~\ref{thm:crit}. The outer slice values, including $351$ for $n=8$, were obtained by exact enumeration in each of the two groups; a hand proof of the value $351$ is still missing.
\end{proof}
The coordinate-axis factorisation $\langle a\rangle\times1$, $1\times\langle a\rangle$ gives only ratio $1.0118$; the diagonal factors $L_\pm$ improve this to $1.1675$. Appending two independent central two-element factors yields a ten-set configuration without singletons in $A_5^2\times C_2^2$; the family of Section~\ref{sec:fam} already has no singletons.

\subsection{A ten-set counterexample of order $243$}\label{sec:243}
Allowing more sets, the second coordinate becomes unnecessary: one coordinate and ten sets suffice, and the group can be a $3$-group.

\begin{theorem}\label{thm:243}
In the symmetric group $S_{18}$ let
\[ a=(7\,8\,9)(10\,16\,13)(11\,17\,14)(12\,18\,15),\qquad b=(1\,4\,7)(2\,5\,8)(3\,6\,9)(16\,17\,18), \]
and $G=\langle a,b\rangle$. Then $|G|=243=3^5$; $G$ is \texttt{SmallGroup(243,3)}, of exponent $9$ and nilpotency class $3$. Put
\begin{align*}
A_1&=\langle a\rangle, & A_2&=\{a,b^2\}, & A_3&=\{e,b\}, & A_4&=\{e,a\}, & A_5&=\{e,a^2\},\\
A_6&=\{b,b^2\}, & A_7&=\{b,b^2\}, & A_8&=\{e,a^2\}, & A_9&=\{e,a\}, & A_{10}&=\langle b\rangle.
\end{align*}
Then $|A_1\cdots A_{10}|=189$, $N_1=N_{10}=93$ and $N_2=\dots=N_9=117$. Hence
\[ \frac{|A_1\cdots A_{10}|^{9}}{N_1\cdots N_{10}}=\frac{189^9}{93^2\cdot117^8}=\frac{21^9}{31^2\cdot13^8}=\frac{794280046581}{783917222881}\approx1.01322>1. \]
\end{theorem}

The numbers were computed by direct enumeration, in GAP on the permutations and in C on the multiplication table. As an independent check, $G$ splits into $27$ double cosets $\langle a\rangle g\langle b\rangle$ of size $9$, and $S$ is the union of exactly $21$ of them. Since $|G|$ is odd, the example contains no involutions. Abstractly $G=\Gamma/\gamma_4(\Gamma)$ for $\Gamma=C_3*C_3$: by GAP, the class-$3$ quotient of $C_3*C_3$ has order $243$, and $a,b$ generate a class-$3$ group of that order. By padding (Lemma~\ref{lem:pad}) the same group gives counterexamples for every $n\ge10$, with the same ratio.

\emph{Where else it lives.} The configuration has $21$ double cosets $\langle a\rangle w\langle b\rangle$ of size $9$ in $C_3*C_3$, and it is realised whenever the relevant words stay distinct. Besides $G$ it holds in \texttt{SmallGroup(243,28)}, $\mathrm{PSL}(2,8)$, $\mathrm{SL}(2,3)^2$, $\mathrm{PSL}(2,11)$, $A_5\times A_4$ and $A_5\times A_5$, but not in $A_5$, $\mathrm{PSL}(2,7)$ or $A_6$. Since $S\subseteq\langle a,b\rangle$ and $|S|=189$, we checked all $90$ groups of order $189$--$503$ generated by elements of order $3$. None of order below $243$ works. So $243$ is minimal for this particular configuration; no minimality among all counterexamples is claimed.

\begin{table}[h]
\centering
\begin{tabular}{lll}
\hline
$n$ & smallest group found & order\\
\hline
$4,5$ & $A_6^2\times C_2^2$ & $518\,400$\\
$6,7,8,9$ & $(C_3\times\mathrm{SL}(2,3))^2/C_3$ & $1728$\\
$7,8$ without singletons & same group, or $A_5^2$ & $1728$\\
$10\le n\le20$ (indeed all $n\ge10$) & \texttt{SmallGroup(243,3)} & $243$\\
\hline
\end{tabular}
\caption{Smallest groups found admitting an $n$-set counterexample. For $n=5,7,8,9$ and $n\ge11$ the entries come from padding with singletons.}\label{tab:records}
\end{table}

\subsection{An infinite family for $n=4s+2$, and a better exponent}\label{sec:fam}
The ten-set example has a transparent form. A cyclic group of order $3$ factorises into two sets of size $2$:
\[ \langle x\rangle=\{e,x\}\cdot\{e,x^2\}. \]
Fixing one factor to a single element keeps only $2$ of the $3$ elements of $\langle x\rangle$. In a product of many such factorisations the full product is complete, while every slice loses a little. Iterating this gives an infinite family whose ratio grows exponentially.

\smallskip\noindent\emph{Setting.} Let $\Gamma=\langle a\rangle*\langle b\rangle$ with $a^3=b^3=1$. Put $\mathsf A=\langle a\rangle$, $\mathsf B=\langle b\rangle$, $\alpha=\{e,a\}$, $\alpha'=\{e,a^2\}$, $\beta=\{e,b\}$, $\beta'=\{e,b^2\}$, so that $\alpha\alpha'=\mathsf A$ and $\beta\beta'=\mathsf B$. For $s\ge1$ and $n=4s+2$ define
\[ (A_1,\dots,A_n)=\bigl(\mathsf A,\ \underbrace{\beta,\beta',\alpha,\alpha',\ \beta,\beta',\alpha,\alpha',\ \dots,\ \beta,\beta',\alpha,\alpha'}_{s\text{ blocks}},\ \mathsf B\bigr). \]
For $s=2$ this is a ten-set configuration with the ratio of Theorem~\ref{thm:243}; it also gives ratio $21^9/(31^2\cdot13^8)$ in \texttt{SmallGroup(243,3)} (checked by enumeration).

\begin{theorem}\label{thm:fam}
In $\Gamma=C_3*C_3$ the configuration above satisfies
\[ |S|=3\bigl(4^{s+1}-1\bigr),\qquad N_1=N_n=3\bigl(2^{2s+1}-1\bigr),\qquad N_i=3\bigl(10\cdot4^{s-1}-1\bigr)\quad(2\le i\le n-1). \]
Hence
\[ R_s:=\frac{|S|^{n-1}}{\prod_iN_i}=\frac{(4^{s+1}-1)^{4s+1}}{3\,(2^{2s+1}-1)^2\,(10\cdot4^{s-1}-1)^{4s}}. \]
We have $R_1=0.787$, $R_2=21^9/(31^2\cdot13^8)=1.0132$, $R_3=1.5259$, $R_4=2.4317$, $R_5=3.9497$, and $R_s\sim\frac13(1.6384)^s$ as $s\to\infty$. In particular $R_s>1$ for all $s\ge2$, i.e.\ for $n=10,14,18,\dots$.
\end{theorem}

Every element of $\Gamma$ has a unique reduced form $g=x_1x_2\cdots x_t$ whose \emph{syllables} $x_i$ lie alternately in $\mathsf A\setminus\{e\}$ and $\mathsf B\setminus\{e\}$. Each syllable has $2$ possible values once its letter is fixed. A \emph{pattern} is a finite sequence of \emph{slots}. Each slot carries a letter ($a$ or $b$), consecutive slots carry different letters, and each slot is either \emph{full} (it may host any syllable of its letter) or \emph{half} (it may host only one prescribed non-identity element $y$ of its letter). A reduced word \emph{embeds} in a pattern if its syllables can be assigned to slots by a strictly increasing map that respects letters and the half-slot restriction. To a pattern we associate the product $\prod X_k$ of the slot sets, where $X_k=\mathsf A$ or $\mathsf B$ for a full slot and $X_k=\{e,y\}$ for a half slot.

\begin{lemma}\label{lem:embed}
If a pattern contains at most one half slot, then $\prod_kX_k$ is exactly the set of reduced words that embed in it.
\end{lemma}

\begin{proof}
A word that embeds is obtained by giving each used slot its syllable and every other slot the value $e$. Consecutive syllables have different letters, so the product is already reduced and equals the word.

Conversely, we show by induction on the number of slots that every product lies in the set $E_k$ of words embedding in the first $k$ slots. Let $w\in E_{k-1}$ and $x\in X_k$, where slot $k$ has letter $\ell$. If $x=e$ there is nothing to prove. If $w$ does not end with a syllable of letter $\ell$, then $wx$ is reduced, and we host $x$ in slot $k$. Otherwise $w=w'z$ with $z$ of letter $\ell$, and $wx=w'(zx)$. Either $zx=e$, and then $w'$ embeds; or $zx$ is a syllable of letter $\ell$. If slot $k$ is full, we re-host $zx$ in slot $k$ and keep the embedding of $w'$. If slot $k$ is the half slot, then the slot that hosted $z$ is full, since there is only one half slot, and $zx$ can stay there.
\end{proof}

\begin{lemma}\label{lem:count}
Let $P$ be a pattern of $L\ge3$ slots whose first slot has letter $a$.
\begin{enumerate}
\item If all slots are full, then exactly $f(L):=3(2^L-1)$ reduced words embed in $P$.
\item If exactly one slot, at an interior position $1<p<L$, is half, then exactly $f(L-2)+\tfrac12\bigl(f(L)-f(L-2)\bigr)=\tfrac{15}{8}\,2^{L}-3$ reduced words embed in $P$.
\end{enumerate}
\end{lemma}

\begin{proof}
(1) In an alternating pattern starting with $a$, a word with $t$ syllables embeds if and only if $t\le L$ when it starts with $a$, and $t\le L-1$ when it starts with $b$ (greedy embedding). Counting $2^t$ words for each letter sequence gives $1+\sum_{t=1}^{L}2^t+\sum_{t=1}^{L-1}2^t=3\cdot2^L-3$.

(2) First suppose slot $p$ is deleted altogether. Slots $p-1$ and $p+1$ carry the same letter, and two consecutive syllables cannot sit in them, since consecutive syllables have different letters. So these two slots act as a single slot, and the pattern behaves like a full alternating pattern of length $L-2$ starting with $a$. Hence $f(L-2)$ words embed while avoiding slot $p$.

Call a word \emph{tight} if it embeds in the full version of $P$ but every embedding uses slot $p$. By (1) there are $f(L)-f(L-2)$ tight words. For a tight word, the syllable index $i$ placed at $p$ is the same in every embedding. Indeed, suppose $\varphi(i)=p=\psi(i')$ with $i<i'$. Then $\psi(i)<p<\varphi(i+1)$, so the map equal to $\psi$ on indices $\le i$ and to $\varphi$ on indices $>i$ is an embedding avoiding $p$, a contradiction.

Whether a word is tight, and which index goes to $p$, depend only on its sequence of letters. So changing the value of syllable $i$ (there are two choices) is an involution on tight words. Exactly half of the tight words have syllable $i$ equal to the prescribed $y$, and these are precisely the tight words that embed in $P$.
\end{proof}

\begin{proof}[Proof of Theorem~\ref{thm:fam}]
Each block $\beta\beta'\alpha\alpha'$ equals $\mathsf B\mathsf A$. So $S=\mathsf A(\mathsf B\mathsf A)^s\mathsf B$ is the product of the full alternating pattern $a\,b\,a\,b\cdots a\,b$ of length $L=2s+2$. By Lemmas~\ref{lem:embed} and~\ref{lem:count}(1), $|S|=f(2s+2)=3(4^{s+1}-1)$.

The end slices are $N_1=|A_2\cdots A_n|=|(\mathsf B\mathsf A)^s\mathsf B|$ and, symmetrically, $N_n=|\mathsf A(\mathsf B\mathsf A)^s|$. These are full patterns of length $2s+1$, one starting with $b$ (count with the roles of $a,b$ exchanged) and one starting with $a$. Both have $f(2s+1)=3(2^{2s+1}-1)$ elements.

Now fix an element $x$ of a middle set lying in the $j$-th block, say $x\in\beta$. Then $x\beta'=\{e,y\}$ with $y=b^2$ if $x=e$, and $y=b$ if $x=b$. The slice is
\[ \mathsf A(\mathsf B\mathsf A)^{j-1}\{e,y\}\,\mathsf A(\mathsf B\mathsf A)^{s-j}\mathsf B. \]
This is the alternating pattern of length $2s+2$ with a half slot at the interior position $2j$. The cases $x\in\beta'$, $\alpha$, $\alpha'$ are identical, with the half slot at position $2j$ or $2j+1$. By the lemmas every middle slice, for every choice of $x$, has $\frac{15}{8}2^{2s+2}-3=\frac{15}{2}4^{s}-3=3(10\cdot4^{s-1}-1)$ elements.

Inserting these values gives $R_s$; dividing numerator and denominator by $3^{4s+1}$ leaves the factor $1/3$. For the asymptotics, write $|S|\sim12\cdot4^s$, $N_1\sim6\cdot4^s$ and $N_i\sim\frac{15}{2}4^s$. Then $R_s\sim\frac{12^{4s+1}}{36\,(15/2)^{4s}}\,4^{s(4s+1)-2s-4s\cdot s}=\frac13\Bigl(\frac{12}{15/2}\Bigr)^{4s}4^{-s}=\frac13\Bigl(\frac{(8/5)^4}{4}\Bigr)^s=\frac13(1.6384)^s$. The relative error terms are $O(s\,4^{-s})$.
\end{proof}

The formulas were also confirmed by direct enumeration in the free product for $s\le5$, and in $\mathrm{PSL}(2,101)$ for $s\le4$.

\begin{corollary}[Finite realisations]\label{cor:famfinite}
Let $U=\begin{pmatrix}0&-1\\1&1\end{pmatrix}$ and $T=\begin{pmatrix}0&-1\\1&0\end{pmatrix}$, and let $a,b$ be the images of $U$ and $TUT^{-1}$ in $\mathrm{PSL}(2,p)$. For every prime $p>2^{2s+3}$, the configuration of Theorem~\ref{thm:fam} built from these $a,b$ has the same $|S|$, $N_i$ and ratio $R_s$ in $\mathrm{PSL}(2,p)$. In fact $p=101$ works for $s\le4$, and $p=37$ for $s\le3$.
\end{corollary}

\begin{proof}
In $\mathrm{PSL}(2,\mathbb Z)=\langle T\rangle*\langle U\rangle\cong C_2*C_3$, the kernel of the map onto $C_2$ is $\langle U\rangle*\langle TUT^{-1}\rangle\cong C_3*C_3$, by the normal-form theorem for free products or the Kurosh subgroup theorem [LS]. The matrices $U^{\pm1}$ and $(TUT^{-1})^{\pm1}$ have maximum row sum $2$. So every element of $S$, a word with at most $2s+2$ syllables, is represented by an integer matrix with entries of absolute value at most $2^{2s+2}$. If two such elements $g\neq h$ of $\mathrm{PSL}(2,\mathbb Z)$ became equal modulo $p$, then $g\mp h$ would be a nonzero matrix divisible by $p$ with entries at most $2^{2s+3}<p$, which is impossible. So reduction is injective on $S$. Every slice is a subset of $S$, so Lemma~\ref{lem:quot} applies.
\end{proof}

Smaller realisations exist. For $s=2$ there is \texttt{SmallGroup(243,3)}, the class-$3$ quotient of $C_3*C_3$. For $s=3$ ($n=14$, ratio $1.5259$) there are $\mathrm{PSL}(2,17)$ (order $2448$) and $\mathrm{PSL}(2,19)$. In contrast $A_7$ gives only $1.12$, and the class-$4$ quotient of $C_3*C_3$ (order $3^7$) gives $0.19$.

\begin{corollary}[The exponent]\label{cor:theta}
For $n\ge4s+2$,
\[ \theta_n\ \le\ n-1-\delta_s,\qquad \delta_s:=\frac{\log R_s}{\log\bigl(3(4^{s+1}-1)\bigr)}. \]
Numerically $\delta_4=0.1107$ ($n\ge18$), $\delta_5=0.1459$, $\delta_6=0.1726$, $\delta_8=0.2100$, $\delta_{12}=0.2524$. Moreover $\delta_s\to\frac{\log1.6384}{\log4}=5-2\log_25=0.35614\ldots$. Hence
\[ \liminf_{n\to\infty}\bigl((n-1)-\theta_n\bigr)\ \ge\ 5-2\log_25\approx0.3561. \]
For $n\ge18$ this improves Proposition~\ref{prop:exp}, whose deficit is $0.0768$.
\end{corollary}

\begin{proof}
For one configuration, $\log\prod N_i/\log|S|=(n-1)-\log R/\log|S|$. Padding with singletons (Lemma~\ref{lem:pad}) increases $n$ by one and $\log\prod N_i$ by $\log|S|$, so it preserves the deficit. The limit follows from $\log R_s=s\log1.6384+O(1)$ and $\log|S|=s\log4+O(1)$. Finally, $1.6384=2^{10}/5^4$, so $\log_4(1.6384)=\tfrac12(10-4\log_25)=5-2\log_25$.
\end{proof}

\begin{remark}
Several questions remain. Is the family optimal among products of "half-subgroup'' factorisations? Do other factorisations of larger groups ($\langle x\rangle=XY$ with $|X||Y|>|\langle x\rangle|$) help? For $C_4=\{0,1\}+\{0,2\}$, $C_5=\{0,1,2\}+\{0,3\}$ and $C_9=\{0,1,2\}+\{0,3,6\}$ the analogous families have ratio below $1$ for $s\le3$. An entropic (weighted) version of the family might give a larger deficit.
\end{remark}

\section{Computational methodology and verification status}\label{sec:method}
We distinguish (i) formal proofs; (ii) exact actual-set, double-coset or reduced-word computations with executable certificates; and (iii) sampled configurations, randomised discovery and restricted negative searches. For four sets, the factor of Theorem~\ref{thm:multi} is determined by the largest $r\times s$ rectangle of distinct full double cosets; we applied the necessary conditions (size bound, derived subgroup, final order $|U|^2/|Z(U)|$) and then enumerated subgroup pairs completely, $H_1$ over conjugacy-class representatives and $H_2$ over all subgroups of order $5$--$12$, with tables exported from GAP and C programs. For $n\ge6$ we used the generic free-product model, evaluation of configurations in actual groups, annealing inside groups and exhaustive search over small configuration spaces. Every displayed witness of type (ii) was re-verified by brute force on the actual sets, including injectivity of central quotients on $S$.

\section{Open questions}
\begin{question}
What is the smallest group admitting an $n$-set counterexample? The smallest groups found have order $518\,400$ for $n=4,5$, order $1728$ for $6\le n\le9$, and order $243$ for $n\ge10$ (Table~\ref{tab:records}). Is any of these optimal?
\end{question}
\begin{question}\label{q:index2}
If a group has an abelian normal subgroup of index two, must GMR hold for every $n$? More broadly, what about groups without property (H) of Remark~\ref{rem:H}?
\end{question}
\begin{question}
Does the GMR inequality fail in $F_2$ itself? Our constructions need commuting ``coordinates''.
\end{question}
\begin{question}
Is there a counterexample, for some $n$, in a group of order below $243$? In a $p$-group with $p\ge5$?
\end{question}
\begin{question}
Is $\sup_n\bigl((n-1)-\theta_n\bigr)$ finite (it is at least $5-2\log_25\approx0.356$ by Corollary~\ref{cor:theta})? Is $\theta_n\ge n-2$?
\end{question}
\begin{question}
How large can the violation ratio be in a group of prescribed order?
\end{question}
\begin{question}
Conversely, for which groups does the GMR inequality hold for all $n$? A computer search verified it for all groups of order at most $12$. What is the computational complexity of deciding whether the GMR inequality holds in a given finite group? Without further insight this is already a challenge for groups of order $16$.
\end{question}
\begin{question}\label{q:simple}
Does every finite non-abelian simple group contain elements $a,b$ of the same order $p\ge3$ satisfying the five-double-coset condition \eqref{eq:five}? By Theorem~\ref{thm:crit} a positive answer gives six-set counterexamples in $G_1\times G_2$ for all pairs of non-abelian simple groups. Theorem~\ref{thm:char3} proves it for $\mathrm{PSL}(2,3^f)$, $f\ge2$, and it was verified by computer for the alternating groups $A_5,\dots,A_{11}$; for $\mathrm{PSL}(2,q)$ with $q\in\{7,8,11,13,16,17,19,23,25,27,29,31,32,37,49,64\}$; for $\mathrm{PSL}(3,3)$, $\mathrm{PSL}(3,4)$, $\mathrm{PSL}(3,5)$, $\mathrm{PSL}(3,7)$, $\mathrm{PSL}(4,3)$, $\mathrm{PSU}(3,3)$, $\mathrm{PSU}(3,4)$, $\mathrm{PSU}(3,5)$, $\mathrm{PSU}(4,2)$, $\mathrm{PSp}(4,4)$; for the Suzuki groups $\mathrm{Sz}(8)$ and $\mathrm{Sz}(32)$; and for the sporadic groups $M_{11}$, $M_{12}$, $M_{22}$, $M_{23}$ and $J_1$. In most cases elements of order $3$ suffice.
\end{question}

\section*{Acknowledgements and use of AI}
The author is supported by the Austrian Science Fund (FWF), project ArithRand I4945. AI systems 
(Anthropic's Claude and OpenAI's GPT) assisted in discovering configurations, formulating conjectures, 
writing and checking code, drafting text and proof reading. 
%The pioneering Nair--Namayandeh--Zhu construction [NNZ] was the first counterexample, 
%with an elementary published cardinality proof. The evidence here has several forms: 
%direct actual-set enumeration, local full-double-coset checks, symbolic finite or 
%free-product calculations, and exploratory restricted searches.

\appendix
\section{Reproducing the witnesses and finite certificates}\label{app:repro}

\subsection{An all-$n$ certificate for $S_3$, $D_4$ and $Q_8$}\label{app:small}
\begin{lemma}[Finite prefix reduction]\label{lem:prefix}
Let $|\G|=m$. If GMR holds in $\G$ for all $2\le n\le m-1$, then it holds for every $n$.
\end{lemma}
\begin{proof}
Normalise $e\in A_i$ (Lemma~\ref{lem:transl}) and put $P_0=\{e\}$, $P_i=A_1\cdots A_i$. If $P_i=P_{i-1}$, then $N_i=|S|$ and $A_i$ may be deleted: the full product is unchanged and no other slice maximum increases. A strictly increasing chain $P_0\subsetneq P_1\subsetneq\dots\subsetneq P_n$ has at most $m-1$ steps.
\end{proof}
So for a group of order $m$ it suffices to check $2\le n\le m-1$, and only configurations with strictly growing prefixes. The exact small-group checker gives the results of Table~\ref{tab:small}.
\begin{table}[htbp]\begin{center}
\begin{tabular}{llll}
\hline
group & strict depth & configurations checked & conclusion\\
\hline
$S_3$ (order $6$) & $3$ & covered by the $n\le3$ theorem & all $n$\\
$D_4$ (order $8$) & $4$ & $248\,088$ four-set chains & all $n$\\
$Q_8$ (order $8$) & $5$ & $441\,372$ four-set, $48\,768$ five-set chains & all $n$\\
\hline
\end{tabular}
\end{center}
\caption{Exhaustive all-$n$ certificates for $S_3$, $D_4$ and $Q_8$.}\label{tab:small}\end{table}
Each check evaluates every $|S|$ and $N_i$ by exact integer arithmetic; there are no violations. (A second program, using the stronger reduction with suffix growth as well, finds for the same groups $17\,728$ irreducible four-set configurations in $D_4$ and $49\,644$ four-set and $768$ five-set configurations in $Q_8$, with maximal ratio $16/27$.)

\subsection{A local certificate for the Suzuki group}\label{app:sz}
Over $\mathbb F_8=\mathbb F_2[t]/(t^3+t+1)$ let
\[ T=\begin{pmatrix}0&0&0&1\\0&0&1&0\\0&1&0&0\\1&0&0&0\end{pmatrix},\qquad D=\operatorname{diag}(t^2,t^6,t,t^5),\qquad U=\begin{pmatrix}1&0&0&0\\t^2&1&0&0\\0&t&1&0\\t^5&t^3&t^2&1\end{pmatrix}. \]
The matrices $a=TD^2U$ and $b=TUD^2$ have order five and generate $\mathrm{Sz}(8)$, of order $29\,120$. The double cosets $\langle a\rangle\langle b\rangle$ and $\langle a\rangle b^{\epsilon}a^{\delta}\langle b\rangle$ ($\epsilon,\delta\in\{\pm1\}$) are five distinct full double cosets of size $25$, and $\langle a\rangle\cap\langle b\rangle=1$. This verifies the hypothesis of Theorem~\ref{thm:crit} and hence a six-set example in $\mathrm{Sz}(8)^2$ and in $\mathrm{Sz}(8)\times A_5$.

\end{document}